\documentclass[10pt]{amsart}
\usepackage{amsmath,amscd,latexsym,verbatim,amssymb}
\usepackage{times}       
 \newcommand{\resp}{{\it resp.} }
\newcommand{\cf}{{\it cf.} }
\newcommand{\ie}{{\it i.e.} }
\newcommand{\eg}{{\it e.g.} }

\newcommand{\N}{\mathbb{N}}
\newcommand{\Q}{\mathbb{Q}}

  \newcommand{\sA}{{\mathcal{A}}}
\newcommand{\sB}{{\mathcal{B}}}
\newcommand{\sC}{{\mathcal{C}}}

\newcommand{\sL}{{\mathcal{L}}}

\newcommand{\sK}{{\mathcal{K}}}

\newcommand{\sO}{{\mathcal{O}}}

\newcommand{\inj}{\hookrightarrow}

 \newcommand{\e}{\frac{1}{p^\infty}}
\newcommand{\f}{-\frac{1}{p^\infty}}
\newcommand{\s}{\sm o}

\font\sm=cmr10 at 10pt

 \newcommand{\colim}{\operatorname{colim}}

\newcommand{\Spec}{\operatorname{Spec}}

\newcounter{spec}

\swapnumbers

\newtheorem{thm}{Theorem}[subsection]
\newtheorem{lemma}[thm]{Lemma}

\newtheorem{cor}[thm]{Corollary}

\theoremstyle{definition}

\numberwithin{equation}{section}

\font\sm=cmr10 at10pt
 at12pt

\renewcommand{\qed}{\hfill $\square$\medskip}
  
\begin{document}

\title[Perfectoid spaces and the homological conjectures.]{Perfectoid spaces and the homological conjectures}
 
\author{Y. Andr\'e}
 
  \address{Institut de Math\'ematiques de Jussieu, 4 place Jussieu, 75005 Paris France.}
\email{yves.andre@imj-prg.fr}
\keywords{homological conjectures, big Cohen-Macaulay algebra, perfectoid algebra, purity}
 \subjclass{13D22, 13H05, 14G20}
  \begin{abstract}  This is a survey of recent advances in commutative algebra, especially in mixed characteristic, obtained by using the theory of perfectoid spaces. An explanation of these techniques and a short account of the author's proof of the direct summand conjecture are included. One then portrays the progresses made with these (and related) techniques on the so-called homological conjectures.
  \end{abstract}
\maketitle

 \begin{sloppypar}

  \section{The direct summand conjecture}
  
\subsection{}  Let $R$ be a Noetherian (commutative) ring and $S$ a finite ring extension, and let us consider the exact sequence of finitely generated $R$-modules
   \begin{equation}\label{1.1} 0 \to R \to S \to S/R\to 0 .\end{equation}    
   
   When does this sequence split? Equivalently, when is $R\to S$ {\it pure}, \ie remains injective after any base change? 
        This holds for instance when $R\to S$ is flat, or when $R$ is a normal $\Q$-algebra, but not in general (the embedding of $\Q[x,y]/(xy)$ in its normal closure gives a counter-example, since it is no longer an embedding modulo $x+y$). 
   
  The direct summand conjecture, formulated by M. Hochster around 1969, claims that {\eqref{1.1} splits whenever $R$ is regular}. 
    Hochster proved it when $R$ contains a field \cite{H1}. R. Heitmann proved it in dimension $\leq 3$ \cite{He}. 
   
  \subsection{}    Recently, the author proved it in general \cite{A2}:
   
   \begin{thm}\label{T1}\eqref{1.1} splits if $R$ is regular.
   \end{thm} 
   
   This has many (non trivially) equivalent forms. One of them is that {\it every ideal of a regular ring $R$ is contracted from every finite (or  integral) extension of $R$}. Another (more indirect) equivalent form is the following statement, which settles a question raised by L. Gruson and M. Raynaud \cite[1.4.3]{RG}\footnote{\cf \cite{O} for the equivalence. Gruson and Raynaud settled the case of a finite extension and outlined that the transition to integral extensions is not a routine exercise.}:
   
     \begin{thm} Any integral extension of a Noetherian ring descends flatness of modules.
   \end{thm}

       We will see one more equivalent form below in the framework of the so-called homological conjectures (\ref{ic}).

   \section{The role of perfectoid spaces}
 
 \subsection{Some heuristics} After Hochster's work \cite{H2}, it is enough to prove the direct summand conjecture in the case when $R$ is a complete unramified local regular ring of mixed characteristic $(0,p)$ and perfect residue field $k$. By Cohen's structure theorem, one may thus assume $R= W(k)[[x_1, \ldots, x_n]]$. 
 
 In characteristic $p$, all proofs of the direct summand conjecture use the Frobenius endomorphism $F$ in some way. In mixed characteristic, $R= W(k)[[x_1, \ldots, x_n]]$ carries a Frobenius-like endomorphism (acting as the canonical automorphism of $W(k)$ and sending $x_i$ to $x_i^p$), which  however does not extend to general finite extensions $S$ of $R$. To remedy this, $p$-adic Hodge theory suggests to ``ramify deeply", by adjoining iterated $p^{th}$ roots of $p, x_1, \ldots, x_n$. Doing this, one leaves the familiar shore of Noetherian commutative algebra for perfectoid geometry, recently introduced by P. Scholze \cite{S1}. 
 
 To begin with, $W(k)$ is replaced by the non-Noetherian complete valuation ring $\sK^{\s} := \widehat{W(k)[p^{\e}]}$. The valuation ring $\sL^{\s}$ of any finite extension $\sL$ of the field $\sK[\frac{1}{p}]$ satisfies:
 \begin{equation}\label{2.1}{ F : \sL^{\s}/p \stackrel{x\mapsto x^p}{\to} \sL^{\s}/p \;\; {\text{\it is surjective, and}} \; \sL^{\s} \; {\text{\it is }} p^{\e}{\text{\it-almost finite etale over}}\;\sK^{\s},} \end{equation}     
this being understood in the context of almost ring theory, introduced by G. Faltings and developped by O. Gabber and L. Ramero \cite{GR}, which gives precise meaning to ``up to $p^{\e}$-torsion"; for instance, $p^{\e}$-almost etaleness means that $p^{\e}\Omega_{\sL^{\s}/\sK^{\s}}=0$. 
 Actually, \cite{GR} is much more general: it deals with modules over a commutative ring up to $\frak k$-torsion, for some idempotent ideal $\frak k$. Going beyond the case of a valuation ideal $\frak k$ will be crucial: beside ``$p^{\e}$-almost" modules, we will have to consider ``$(pg)^{\e}$-almost" modules for some ``geometric" discriminant $g$.  

   \subsection{Perfectoid notions} In perfectoid geometry, one works with certain Banach\footnote{here and in the sequel, one can work with any perfectoid field $\sK$ (of mixed char. $(0,p)$), -  \ie  complete, non-discretely  valued, and such that $F$ is surjective on $\sK^{\s}/p$. An extensive dictionary between the langage of commutative algebra and the language of non-archimedean functional analysis is presented in \cite[2.3.1]{A1}.} $\sK$-algebras $\sA$. One denotes by $\sA^{\s}$ the $\sK^{\s}$-subalgebra of power-bounded elements. One says that $\sA$ is {\it uniform} if $\sA^{\s}$ is bounded, and that $\sA$ is {\it perfectoid} if it is uniform and $F : \sA^{\s}/p \stackrel{x\mapsto x^p}{\to} \sA^{\s}/p$ is surjective.
 An example which plays a crucial role in the sequel is $\hat\cup_i W(k)[p^{\frac{1}{p^i}}][[x_1^{\frac{1}{p^i}} , \ldots, x_n^{\frac{1}{p^i}}]][\frac{1}{p}]\,$, a deeply ramified avatar of $R$. 
 Morphisms of perfectoid algebras $\sA \to \sB$ are continuous algebra homomorphisms (one then says that $\sB$ is a perfectoid $\sA$-algebra).  
 
 Perfectoid algebras enjoy three fundamental stability properties \cite{S1}:

\subsubsection{\emph{Tensor product}}\label{2.2.1} {\it If $\sB$ and $\sC$ are perfectoid $\sA$-algebras, so is $\sB\hat\otimes_\sA \sC$, and $(\sB\hat\otimes_\sA \sC)^{\s}$ is $p^{\e}$-almost isomorphic to $ \sB^{\s}\hat\otimes_{\sA^{\s}} \sC^{\s}$.}

\subsubsection{\emph{Localization}}\label{2.2.2} {\it The ring of functions $\sA\{\frac{f}{g}\}$ on the  subset of the perfectoid space ${\rm{Spa}}(\sA, \sA^{\s})$ where $\vert  f\vert \leq \vert g\vert $ holds is perfectoid, and $\sA\{\frac{f}{g}\}^{\s}$ is $p^{\e}$-almost isomorphic to $\sA^{\s}\langle (\frac{f'}{g'})^{\e}\rangle$ for some approximations $f', g'$ of $f, g$ which admit iterated $p^{th}$-roots in $\sA$.}

\subsubsection{\emph{Finite etale extension}}\label{2.2.3} {\it Any finite etale extension $\sB$ of $\sA$ is perfectoid, and $\sB^{\s}$ is a $p^{\e}$-almost finite etale extension of $\sA^{\s}$.}

 This generalization of \ref{2.1} to perfectoid algebras is Faltings's ``almost purity theorem" \cite{Fa} as revisited by Scholze \cite{S1} and Kedlaya-Liu \cite{KL}.  
 
 \medskip {\small Let us explain how the second assertion of \ref{2.2.3} follows from  the first following \cite[3.4.2]{A1}.  The idea is to reduce to the case when $\sB $ is Galois over $\sA $ with Galois group $G$, \ie $\sB^G = \sA $ and $\sB\otimes_{\sA} \sB\stackrel{\sim}{\to} \prod_G  \sB$. This implies $\sB^{\s G} = \sA^{\s}$. On the other hand,
 since $\sB $ is a finitely generated projective ${\sA }$-module, $\sB  \otimes_{\sA } \sB  =  \sB  \hat\otimes_{\sA } \sB $, and one deduces from \ref{2.2.1} (assuming $\sB$ perfectoid) that $\sB^{\s}\hat\otimes_{\sA^{\s}} \sB^{ \s} \to \prod_G  \sB^{\s}$ is a $p^{\e}$-almost isomorphism. To get rid of the completion, one passes modulo $p^m$:  $\sB^{ \s}/p^m$ is almost Galois over $\sA^{ \s}/p^m$, hence almost finite etale, and a variant of Grothendieck's ``\'equivalence remarquable" \cite[5.3.27]{GR} allows to conclude that  $ \sB^{ \s} $ it itself almost finite etale over $ \sA^{ \s} $. }

 \subsection{Direct summand conjecture: the case when $S[\frac{1}{p}]$ is etale over $R[\frac{1}{p}]$.}\label{sc} Let us go back to the direct summand conjecture for $R = W(k)[[x_1, \ldots, x_n]]$. The special case when $S[\frac{1}{p}]$ is etale over $R[\frac{1}{p}]$ was settled  by B. Bhatt \cite{Bh} and by K. Shimomoto \cite{Sh}, using   \ref{2.2.3}. Here is a slightly different account, suitable to the sequel. 
   
 Let us consider  the perfectoid algebra $\sA = \hat\cup_i W(k)[p^{\frac{1}{p^i}}][[x_1^{\frac{1}{p^i}} , \ldots, x_n^{\frac{1}{p^i}}]][\frac{1}{p}]$ and notice that $\sA^{\s} =  \hat\cup_i W(k)[p^{\frac{1}{p^i}}][[x_1^{\frac{1}{p^i}} , \ldots, x_n^{\frac{1}{p^i}}]]$ is a faithfully flat extension of $R$. By assumption $\sB := S\otimes_R \sA$ is finite etale over $\sA$, hence by \ref{2.2.3}, $\sB^{\s}$ is $p^{\e}$-almost finite etale, hence almost pure (in the sense: almost universally injective) over $\sA^{\s}$. 
 A fortiori, \eqref{1.1} almost splits after tensoring with $\sA^{\s}$. In other words, if  $e\in {\rm Ext}^1_R(S/R, R)$ denotes the class corresponding to \eqref{1.1}, then $p^{\e}(e\otimes 1)=0$ in $ {\rm Ext}^1_R(S/R, R)\otimes_R \sA^{\s}  \cong {\rm Ext}^1_{\sA^{\s}}((S\otimes_R  \sA^{\s})/  \sA^{\s}, \sA^{\s})$. One concludes that $e=0$ by the following general elementary lemma (applied to $M= Re$ and $\frak K= p^{\e} \sA^{\s}$):
   
      \begin{lemma}\label{L2}\cite[1.1.2]{A2} Let $R$ be a local Noetherian ring, $M$ a finitely generated  $R$-module, $A$ a faithfully flat $R$-algebra. Let $\frak K $ be an idempotent ideal of $A$ such that $\frak K . M_A = 0$ and $R\cap  \frak K   \neq 0$. Then $M = 0$. 
  \end{lemma}

\subsection{The perfectoid Abhyankar lemma}
 In the general case, $S\otimes_R \sA$ is no longer etale over $\sA$: one must take into account a discriminant $g\in R$ of $S[\frac{1}{p}]$  over $R[\frac{1}{p}]$.
 This suggests to try to generalize \ref{2.2.3} to {\it ramified} extensions of perfectoid algebras. It turns out that this is possible, provided one extracts suitable roots of $g$ in the spirit of Abhyankar's lemma.    
 
This leads to replace everywhere ``$p^{\e}$-almost" by ``$(pg)^{\e}$-almost", thereby extending the basic setting of almost ring theory beyond the usual situation of a non-discrete valuation ring. This also leads to introduce the notion of {\it almost perfectoid algebra}, where 
$F : \sA^{\s}/p \stackrel{x\mapsto x^p}{\to} \sA^{\s}/p$ is only assumed to be $(pg)^{\e}$-almost surjective \cite[3.5.4]{A1}.
  
 \begin{thm}\label{T2}\cite{A1} Let $\sA$ be a perfectoid $\sK$-algebra, which contains a compatible system of $p$-power roots $g^{\frac{1}{p^j}}$ of some non-zero-divisor $g\in \sA^{\s}$. Let $\sB'$ be a finite etale $\sA[\frac{1}{g}]$-algebra. Let $\sB$ be the integral closure of $g^{\f}\sA$ in $\sB'$, so that $\sB[\frac{1}{g}]= \sB'$. 
 
  Then $\sB$ is almost perfectoid, and for any $n$, $\sB^{\s}/p^m$ is $(pg)^{\e}$-almost finite etale (hence $(pg)^{\e}$-almost flat) over $\sA^{\s}/p^m$. 
  \end{thm} 
 
 (If $g=1$, one may use \cite[5.3.27]{GR} again to conclude that  $ \sB^{ \s} $ is itself almost finite etale over $ \sA^{ \s} $ and recover \ref{2.2.3}). 
 
 \smallskip
The basic idea is to look at the pro-system of algebras of functions $\sA^j := \sA\{\frac{p^j}{g}\} $ on complements of tubular neighborhoods of the hypersurface $g=0$ in the perfectoid space ${\rm{Spa}}(\sA, \sA^{\s})$, \resp at the pro-system $\sB^j := \sB'\otimes_{\sA[\frac{1}{g}]}\sA\{\frac{p^j}{g}\} $. Each $\sA^j$ is perfectoid (\ref{2.2.2}), and each $\sB^j$ is finite etale over $\sA^j$, hence perfectoid; moreover $\sB^{j\s}$ is $p^{\e}$-almost finite etale over $\sA^{j\s}$ by almost purity (\ref{2.2.3}). One can show that 
 $\sB^{\s}$ is isomorphic to $\lim \sB^{j\s}$, and that the latter has the asserted properties, 
 
\smallskip However, in the sequel, the identification of $(\lim \sB^{j\s})[\frac{1}{p}]$ with the integral closure of $g^{\f}\sA $ in $\sB'$ plays no role; changing notation, we will set $\sB := (\lim \sB^{j\s})[\frac{1}{p}]$, which is a uniform Banach algebra, and sktech the  {\it proof that $\sB$ is almost perfectoid and that for every $m$, $\sB^{\s}/p^m$ is $(pg)^{\e}$-almost finite etale  over $\sA^{\s}/p^m$}. 
 
 {\small \medskip  The proof involves six steps. 

\smallskip 1) For any  $r\in \N[\frac{1}{p}]$, {\it $\lim (\sA^{j\s}/p^r)$ is $(pg)^{\e}$-almost isomorphic to $\sA^{\s}/p^r$}. This is essentially Scholze's perfectoid version of Riemann's extension theorem \cite[II.3.1]{S2} (hint: if $R$ denotes $\sA^{\s}/p^r$ for short, $\sA^{j\s}/p^r$ is $p^{\e}$-almost isomorphic to $R^j := R[(\frac{p^j}{g})^{\e}]$; the key idea is that for $j'\geq j+ r p^k , R^{j'}\to R^j$ factors through $R^{jk} := \sum_{s\leq \frac{1}{p^k}} R(\frac{p^j}{g})^{s} $, so that $\lim R^j \cong \lim R^{jk}$; on the other hand, the kernel and cokernel of $R \to R^{jk}$ are killed by  $g$ raised to a power which tends to $0$ when $j, k \to \infty$). 
 Passing to the limit $r\to \infty$, it follows that $\lim  \sA^{j\s}=g^{\f}\sA,$ and this also holds under the weaker assumption that $\sA$ is almost perfectoid,  \cf \cite[4.2.2]{A1}.
 
\smallskip 2) {\it $\lim^1 (\sB^{j\s}/p^r)$ is $(pg)^{\e}$-almost zero}. 
The technique is similar to the one in 1), \cf \cite[4.4.1]{A1}. 

\smallskip 3) {\it  $  \lim (\sB^{j\s}/p) \stackrel{F}{\to} \lim (\sB^{j\s}/p)\, $ is $(pg)^{\e}$-almost surjective}. Indeed, taking the limit of the exact sequence 
$ 0\to \sB^{j\s}/p^{\frac{p-1}{p}} \to  \sB^{j\s}/p\to  \sB^{j\s}/p^{\frac{1}{p}} \to 0$, one deduces from $2)$ (for $r= \frac{p-1}{p})$, that $\lim (\sB^{j\s}/p) \to \lim (\sB^{j\s}/p^{\frac{1}{p}})$ is $(pg)^{\e}$-almost surjective; on the other hand, $\sB^{j\s}/p^{\frac{1}{p}} \to \sB^{j\s}/p$ is an isomorphism because $\sB^{j}$ is perfectoid.

\smallskip 4) {\it $\sB$ is almost perfectoid}.  From 3), it suffices to show that the natural map $\sB^{\s}/p =  (\lim \sB^{j\s})/p \to \lim (\sB^{j\s}/p)$ is a $(pg)^{\e}$-almost isomorphism. It is easy to see that it is injective \cite[2.8.1]{A1}, and on the other hand, the composition $     \lim_{F, j} (\sB^{j\s}/p) = \sB^{\flat\s}   \to (\lim \sB^{j\s})/p \to \lim (\sB^{j\s}/p) $ is almost surjective by 3).

\smallskip 5) {\it $\sB \to \sB^{j }$ factors through a $g^{\e}$-almost isomorphism $\sB\{\frac{p^j}{g}\}    \stackrel{a}{\cong} \sB^{j }$ }. The factorization comes from the fact that $ \sB^{j } \cong  \sB^{j }\{\frac{p^j}{g}\} $, and one constructs an almost inverse to $\sB\{\frac{p^j}{g}\}    \stackrel{a}\to  \sB^{j }$ as follows (\cf \cite[4.4.4]{A1}): the integral closure of $\sA^{\s}$ in $\sB'$ maps to the integral closure $\sA^{j\s}$ in $\sB^j$, which is $p^{\e}$-almost $\sB^{j\s}$ (by almost purity). Passing to the limit and inverting $pg$, one gets a morphism $\sB' \stackrel{\delta}{\to} \sB[\frac{1}{g}]$, and the sought for inverse is induced by $\delta\otimes 1_{\sA^j}$.   
 
\smallskip 6)  {\it for every $m$, $\sB^{\s}/p^m$ is $(pg)^{\e}$-almost finite etale  over $\sA^{\s}/p^m$}. This is the decisive step: how to keep track of  the almost finite etaleness of $\sB^{j\s}$ over $\sA^{j\s}$ at the limit? 
  As in \ref{2.2.3}, the idea is to reduce to the case when $\sB'$ is Galois over $\sA[\frac{1}{g}]$ with Galois group $G$, \ie $(\sB')^G = \sA[\frac{1}{g}]$ and $\sB'\otimes_{\sA[\frac{1}{g}]} \sB' \stackrel{\sim}{\to} \prod_G  \sB'$. It follows that $\sB^j$ is $G$-Galois over $\sA^j$, and (as in \ref{2.2.3}) that $\sB^{j \s}\hat\otimes_{\sA^{j \s}} \sB^{j \s} \to \prod_G  \sB^{j \s}$ is a $p^{\e}$-almost isomorphism.   
  Passing to the limit, (and taking into account the facts that, by 4) and 5), $\sB \hat\otimes_{\sA} \sB$ is an almost perfectoid algebra and $\sB \hat\otimes_{\sA} \sB\{\frac{p^j}{g}\} \cong \sB^j \hat\otimes_{\sA^j}  \sB^j $, so that one can apply 1)), one concludes that $\sA^{\s} \to \sB^{\s G}$ and $\sB^{  \s}\hat\otimes_{\sA^{  \s}} \sB^{ \s} \to \prod_G  \sB^{ \s}$ are $(pg)^{\e}$-almost isomorphisms. To get rid of the completion, one passes modulo $p^m$:  $\sB^{ \s}/p^m$ is almost Galois over $\sA^{ \s}/p^m$, hence $(pg)^{\e}$-almost finite etale (in constrast to \ref{2.2.3}, one cannot conclude that $\sB^{ \s} $ is $(pg)^{\e}$-almost finite etale over $\sA^{ \s} $ since $p$ may not belong to $(pg)^{\e}$).} \qed

\subsection{Infinite Kummer extensions}
  
  In order to extend the strategy of \ref{sc} to the general case by means of the perfectoid Abhyankar lemma, one has first to adjoin to the perfectoid algebra $  \hat\cup_i W(k)[p^{\frac{1}{p^i}}][[x_1^{\frac{1}{p^i}} , \ldots, x_n^{\frac{1}{p^i}}]][\frac{1}{p}]$ the iterated $p^{th}$-roots of a discriminant $g$.   
  At finite level $i$, it seems very difficult to control the extension $(W(k)[p^{\frac{1}{p^i}}][[x_1^{\frac{1}{p^i}} , \ldots, x_n^{\frac{1}{p^i}}]][g^{\frac{1}{p^i}},\frac{1}{p}])^{\s}$, which is bigger than $W(k)[p^{\frac{1}{p^i}}][[x_1^{\frac{1}{p^i}} , \ldots, x_n^{\frac{1}{p^i}}]][g^{\frac{1}{p^i}}]$, but things turn easier at the infinite level, thanks to the perfectoid theory.  
  
   \begin{thm}\label{T3}\cite[2.5.2]{A2} Let $\sA$ be a perfectoid $\sK$-algebra, and let $g\in \sA^{\s}$ be a non-zero divisor. Then for any $n$, $\sA\langle g^{\e}\rangle^{\s}/p^m$ is $p^{\e}$-almost faithfully flat over $\sA^{\s}/p^m$.
    \end{thm}
   
   The basic idea is to add one variable $T$, consider the perfectoid algebra $\sC := \sA\langle T^{\e}\rangle$ and look at the ind-system of algebras of functions $\sC_i := \sC\{\frac{T-g}{p^i} \} $ on  tubular neighborhoods  of the hypersurface $T= g$ in the perfectoid space ${\rm{Spa}}(\sC, \sC^{\s})$. 
   
 {\small \medskip  The proof involves three steps. 
   
 \smallskip  7) {\it $\sA\langle g^{\e}\rangle^{\s}$ is $p^{\e}$-almost isomorphic to $\widehat\colim_i \,\sC_i^{\s}$}. This is an easy consequence of the general fact that for any uniform Banach algebra $\sB$ and $f\in \sB^{\s}$, $\widehat\colim_i \sB\{\frac{f}{p^i}\}^{\s}$ is $p^{\e}$-almost isomorphic to $(\sB/f\sB)^{\s}$, \cf \cite[2.9.3]{A1}.
 
  \smallskip  8) {\it $\sC^{\s}$ contains a compatible system of $p$-power roots of some non-zero-divisor $f_i$ such that $\sC_i ^{\s}$  is $p^{\e}$-almost isomorphic to $\widehat\colim_j \, \sC\langle  U\rangle^{\s}/ (p^{\frac{1}{p^j}} U - f_i^{\frac{1}{p^j}}) $.} This is one instance of Scholze's approximation lemma \cite[6.7]{S1}; one may assume $f_i \equiv T-g \,\mod p^{\frac{1}{p}}$.
   
     \smallskip  9) {\it $\sC\langle  U\rangle^{\s}/ (p^{\frac{1}{p^j}} U - f_i^{\frac{1}{p^j}}, p^m) $ is faithfully flat over $\sA^{\s}/p^m$.} One may replace $p^m$ by any positive power of $p$, \eg $p^{\frac{1}{p^{j+1}}}$. Since $\sC$ is perfectoid, there is $g_{ij}\in  \sC^{\e}\rangle^{\s}$ such that $ g_{ij}^{p ^j} \equiv g  \,\mod p$. Then  $f_i^{\frac{1}{p^j}} \equiv T^{\frac{1}{p^j}} -g_{ij}  \,\mod p^{\frac{1}{p^{j+1}}},$ and $\sC\langle  U\rangle^{\s}/ (p^{\frac{1}{p^j}} U - f_i^{\frac{1}{p^j}}, p^{\frac{1}{p^{j+1}}}) \cong (\sA^{\s}/p^{\frac{1}{p^{j+1}}})[T^{\e}, U]/ (T^{\frac{1}{p^j}} -g_{ij})$, a free $\sA^{\s}/p^{\frac{1}{p^{j+1}}}$-module. \qed}

 \subsection{Conclusion of the proof of the direct summand conjecture}\label{conc} One chooses  $g\in R$ such that $S[\frac{1}{pg}]$ is etale over $R[\frac{1}{pg}]$. 
 One then follows the argument of \ref{sc}, replacing   $\sA =  \hat\cup_i W(k)[p^{\frac{1}{p^i}}][[x_1^{\frac{1}{p^i}} , \ldots, x_n^{\frac{1}{p^i}}]][\frac{1}{p}]$  by the ``infinite Kummer extension" $\sA  :=   (\hat\cup_i W(k)[p^{\frac{1}{p^i}}][[x_1^{\frac{1}{p^i}} , \ldots, x_n^{\frac{1}{p^i}}]][\frac{1}{p}])\langle g^{\e}\rangle $. One introduces the finite etale extension $\sB' := S\otimes_R \sA[\frac{1}{g}]$ of $\sA[\frac{1}{g}]$, and one replaces  $\sB = S\otimes_R \sA$ in \ref{sc} by $\sB := (\lim \sB^{j\s})[\frac{1}{p}]$, noting that $\sB^{\s}$ is the $S$-algebra $\lim \sB^{j\s}$. Translating Th. \ref{T2} and Th. \ref{T3} in terms of almost purity modulo $p^m$, and reasoning as in \ref{sc} (using the same lemma), one obtains that \eqref{1.1} splits modulo $p^m$ for any $n$ \cite[\S 3]{A2}. A Mittag-Leffler argument (retractions of $R/p^m\to S/p^m$ form a torsor under an artinian $R/p^m$-module \cite[p. 30]{H1}) shows that \eqref{1.1} itself splits. \qed   
 
  \subsection{Derived version} In \cite{Bh2}, B. Bhatt revisits this proof and proposes a variant, which differs in the analysis of the pro-system $\sA^{j\s}/p^r$ occurring in the proof of Th. \ref{T2}: he strengthens step 1) by showing that the pro-system of kernels and cokernels of $(\sA/p^r)_j \to (\sA^{j\s}/p^r)_j$  is pro-isomorphic to a pro-system of $(pg)^{\e}$-torsion modules; this allows to apply various functors before passing to the limit $j\to \infty$, whence a gain in flexibility. More importantly, he obtains the following derived version of the direct summand conjecture, which had been conjectured by J. de Jong:  
 
\begin{thm}\label{T4}\cite{Bh2} Let $R$ be a regular Noetherian ring, and $f  : X  \to \Spec R$ be a proper surjective morphism. Then the
map $ R \to   R\Gamma(X, \sO_X)$ splits in  the derived category $D(R)$. \end{thm}   
 
 \section{Existence of (big) Cohen-Macaulay algebras}
      
  \subsection{Cohen-Macaulay rings, and Cohen-Macaulay algebras for non-Cohen-Macaulay rings}  Let $S$ be a local Noetherian ring with maximal ideal $\frak m$ and residue field $k$. Recall that a sequence $\underline x = (x_1, \ldots, x_n)$ in $\frak m$ is {\it secant}\footnote{following Bourkaki's terminology (for instance); it is also often called ``part of a system of parameters", although ${\rm{gr}}_{\underline x}S$ may not be a polynomial ring in the ``parameters" $x_i$  (it {\it is} if $\underline x$ is regular).} if $\dim S/(\underline x) = \dim S - n$, and {\it regular} if for every $i$, multiplication by $x_i$ is injective in $S/(x_1, \ldots, x_{i-1})S$.  Any regular sequence is secant, and if the converse holds, $S$ is said to be {\it Cohen-Macaulay}. Regular local rings have a secant sequence generating $\frak m$, and are Cohen-Macaulay. 
      
         Cohen-Macaulay rings 
        form the right setting for Serre duality and the use of local homological methods in algebraic geometry, and have many applications to algebraic combinatorics \cite{BH}. When confronted with a non-Cohen-Macaulay ring $S$, one may try two expedients:
        
    \smallskip    $1)$ {\it Macaulayfication}: construct a proper birational morphism $X \to \Spec S$ such that all local rings of $X$ are Cohen-Macaulay. This weak resolution of singularities, introduced by Faltings, has been established in general by T. Kawasaki \cite{Ka}. Hovewer, secant sequences in $S$ may not remain secant (hence not become regular) in the local rings of $X$; this motivates the second approach:
        
        \smallskip     $2)$ {\it Construction of a Cohen-Macaulay algebra\footnote{since it would be too restrictive to impose that $C$ is Noetherian, one often speaks of ``big" Cohen-Macaulay algebra.}}: an $S$-algebra $C$ such that any secant sequence of $S$ becomes regular in $C$, and $\frak m C\neq C$.      
      
   The existence of Cohen-Macaulay algebras implies the direct summand conjecture: indeed, if $C$ is a Cohen-Macaulay algebra for a finite extension $S$ of a regular local ring $R$, it is also a Cohen-Macaulay $R$-algebra; this implies that $R\to C$ is faithfully flat, hence pure, and so is  $R\to S$.

            \subsection{Constructions of Cohen-Macaulay algebras}   The existence of a (big) Cohen-Macaulay algebra was established by Hochster and C. Huneke under the assumption that $S$ contains a field  \cite{HH}. One may assume that $S$ is a complete local domain. In char. $p>0$, one may then take $C$ to be the {\it total integral closure of $S$} (\ie the integral closure of $S$ in an algebraic closure of its field of fractions). This is no longer true in the case of equal char. $0$, which can nevertheless be treated by reduction to char. $p>>0$ using ultraproduct techniques. 
      
      The remaining case of mixed characteristic was settled in \cite{A2}, using the same perfectoid methods, so that one has: 
      
      \begin{thm}\label{T5} Any local Noetherian ring $S$ admits a (big) Cohen-Macaulay algebra $C$.
      \end{thm} 
      
  {\small In the case of a complete local domain $S$ of char $(0,p)$ and perfect residue field $k$ (to which one reduces), one proceeds as follows. Cohen's theorem allows to present $S$ as a finite extension of $R = W[[x_1, \ldots, x_n]]$. 
      One first considers the $R$-algebra $\sA  :=   (\hat\cup_i W(k)[p^{\frac{1}{p^i}}][[x_1^{\frac{1}{p^i}} , \ldots, x_n^{\frac{1}{p^i}}]][\frac{1}{p}])\langle g^{\e}\rangle^{\s} $ and the $S$-algebra $\sB^{\s} = \lim \sB^{j\s}$ as above \ref{conc}. It follows from Th. \ref{T2} and Th. \ref{T3} that $\sB^{\s}$ is $(pg)^{\e}$-almost isomorphic to a faithfully flat $R$ algebra modulo any power of $p$. From there, one deduces that the sequence $(p, x_1, \ldots, x_n)$ is ``$(pg)^{\e}$-almost regular" in $\sB^{\s}$. 
    
    To get rid of ``almost", Lemma \ref{L2} is no longer sufficient: instead, one uses Hochster's technique of monoidal modifications \cite{H3}\cite{HH}. After $\frak m$-completion, one gets a $S$-algebra $C$ in which $(p, x_1, \ldots, x_n)$, as well as any other secant sequence of $S$, becomes regular. }\qed 
          
   Subsequently, using the tilting equivalence between perfectoid algebras in char. $0$ and in char.$p\,$ and applying Hochster's  modifications in char.$p\,$ rather than in char.$0$, K. Shimomoto \cite{Sh2} shows that {\it in Th. \ref{T5}, in mixed characteristic,  $C$ can be taken to be perfectoid\footnote{actually, Shimomoto's result is slightly weaker, but can be enhanced.}.  In particular, if $S$ is regular, it admits a perfectoid faithfully flat algebra} (one may speculate about the converse\footnote{this is settled in forthcoming work by Bhatt, Ma, Iyangar.}).

  \subsection{Finite and fpqc covers} Since Cohen-Macaulay algebras for regular local rings are faithfully flat, Th. \ref{T5} implies \cite{A2}:
  
  \begin{thm} Any finite cover of a regular scheme is dominated by a faithfully flat quasi-compact cover.
  \end{thm} 
 
 If regularity is omitted, $\Spec(\Q[x,y]/(xy))$ and its normalization provide a counter-example.

   \section{Homological conjectures} 
  
 \subsection{Origins from intersection theory} Under the influence of M. Auslander, D. Buchsbaum and J.-P. Serre, commutative algebra has shifted  in the late 50s from the study of ideals of commutative rings to the homological study of modules (\cf their characterization of regular local rings by the existence of finite free resolutions for any finitely generated module, \resp for the residue field). 
 
 Serre proved that for any three prime ideals $\frak p, \frak q, \frak r$ of a regular local ring $R$ such that $\frak r$ is a minimal prime of $\frak p + \frak q$, ${\rm{ht}}\, \frak r \leq {\rm{ht}}\, \frak p + {\rm{ht}}\, \frak q$ \cite{Se}. The special case $\frak r= \frak m$ can be amplified: for any ideals $I, J$ of $R$ such that $I+J$ is $\frak m$-primary, $\dim R/I+ \dim R/J \leq \dim R$.  
 
 This is no longer true if $R$ is not regular, and attempts to understand the general situation led to the so-called homological conjectures \cf \cite[ch. 9]{BH}, \cite{H4}.

  \subsection{Intersection conjectures}\label{ic} Let $(R, \frak m)$ be local Noetherian ring. 
  
  The first ``intersection conjecture" was proposed by C. Peskine and L. Szpiro \cite{PS}, proved by them when $R$ contains a field by reduction to char. $p$ and Frobenius techniques, and later proved in general by P. Roberts using $K$-theoretic methods \cite{Ro}. It states that {\it if $M, N$ are finitely generated $R$-modules such that $M\otimes N$ has finite length, then $\dim N \leq {\rm{pd}}\, M$.} It implies that {\it $R$ is Cohen-Macaulay if and only if there is an $R$-module $M$ of finite length and finite projective dimension} (in the spirit of  Serre's characterization of regular rings, which is the case $M=k$), resp. {\it if there is an $R$-module $M$ of finite injective dimension}.
  
  {\small Indeed, one may take $N=R$ and deduce that $\dim R\leq  {\rm{pd}}\,M$; by the Auslander-Buchsbaum formula, ${\rm{pd}}\,M  = {\rm{depth}}\, R - {\rm{depth}}\, M$, so that the inequality ${\rm{depth}}\, R \leq \dim R$ is an equality. The second assertion follows from the fact that $ {\rm{id}}\,M  = {\rm{depth}}\, R$ \cite{Ba}.}
  
 \smallskip The ``new intersection conjecture", also proved by Peskine, Szpiro \cite{PS} and Roberts \cite{Ro}, states that {\it for any non exact complex $F_\bullet$ of free $R$-modules concentrated in degrees $[0, s]$ with finite length homology,  $s > \dim R$}.  

  The ``improved new intersection conjecture" is a variant due to E. Evans and P. Griffith \cite{EG}, in which the condition on $F_\bullet$ is ``slightly" relaxed: the $H_{i>0}$ are of finite length and there exists a primitive cyclic submodule of $H_0$ of finite length. 
    They proved it, assuming the existence of (big) Cohen-Macaulay algebras\footnote{if $R$ is Cohen-Macaulay, the Buchsbaum-Eisenbud criterion gives a condition for the exactness of $F_\bullet$ in terms of codimension of Fitting ideals of syzygies. In general, the same condition guarantees that  $F_\bullet \otimes C$ is exact for any Cohen-Macaulay $R$-algebra $C$ \cite[9.1.8]{BH}.}, and showed that it implies  
 their ``syzygy conjecture".  In spite of appearances, the passage from the new intersection conjecture to its ``improved" variant is no small step\footnote{$K$-theoretic techniques failed to make the leap.}: in fact, according to Hochster \cite{H4} and S. Dutta \cite{D}, the latter is {\it equivalent to the direct summand conjecture}. 
 
 On the other hand, in the wake of the new intersection conjecture (and motivated by the McKay correspondence in dimension $3$ and 
the ``fact" that threefold flops induce equivalences of derived categories), 
 T. Bridgeland and S. Iyengar obtained a refinement of Serre's criterion for regular rings assuming the existence of Cohen-Macaulay algebras  \cite[2.4]{BI}.

By Th. \ref{T1} and Th. \ref{T5}, the improved new intersection conjecture and the Bridgeland-Iyengar criterion thus hold inconditionally:  

\begin{thm} Let $R$ be a Noetherian local ring and $F_\bullet$ be a complex of finitely generated free $R$-modules concentrated in degree $[0, s]$, such that $H_{>0}(F_\bullet)$ has finite length.    
 \begin{enumerate} 
 \item If $H_0(F_\bullet)$ contains a cyclic $R$-submodule of finite length not contained in $\frak m H_0(F_\bullet)$, then $s\geq \dim R$.
    \item If $H_0(F_\bullet)$ has finite length and contains $k$ as a direct summand, and $s = \dim R$, then $R$ is regular. 
 \end{enumerate}
\end{thm} 
 
And so does the syzygy conjecture: 

\begin{thm} Let $R$ be a Noetherian local ring and $M$ a finitely generated $R$-module of finite projective dimension $s$. Then for $i\in \{1, \ldots, s-1\}$, the $i$-th syzygy module of $M$ has rank $\geq i$.
\end{thm}

   \subsection{Further work around the homological conjectures using perfectoid spaces} 
 
 \subsubsection{} In \cite{HeMa1}, Heitmann and L. Ma show that Cohen-Macaulay algebras can be constructed in a way compatible with quotients $S\to S/\frak p$ by primes of height one. Using arguments similar to Bhatt's derived techniques, they deduce the vanishing conjecture for maps of Tor \cite{HH}: 
 
 \begin{thm}\cite{HeMa1} Let $R\to S\to T$ be morphisms such that the composed map is a local morphism of mixed characteristic regular local rings, and $S$ is a finite
torsion-free extension of $R$. Then for every $R$-module  $M$ and every $i$, the map ${\rm{Tor}}^R_i(M, S) \to {\rm{Tor}}^R_i(M, T)$ vanishes. 
 \end{thm} 
 
 They obtain the following corollary, which generalizes results by Hochster and J. Roberts, J.-F. Boutot et al.:
 
  \begin{cor}\cite{HeMa1}  Let $R\inj S $ be a pure, local morphism, with $S$ regular. Then $R$ is pseudo-rational, hence Cohen-Macaulay.  
 \end{cor}

 \subsubsection{} In \cite{MaS}, Ma and K. Schwede define and study perfectoid multiplier/test ideals in mixed characteristic, and use them to bound symbolic powers of ideals in regular domains in terms of ordinary powers:   
 
 \begin{thm}\cite{MaS} Let $R$ be a regular excellent Noetherian domain and let $I\subset R$ be
a radical ideal such that each minimal prime of $I$ has height $\leq h$. Then for every $n$, 
 $I^{(hn)} \subset I^n$.\end{thm} 
    Here  $I^{(hn)} $ denotes the ideal of elements of $R$ which vanish generically to order $hn$ at $I$. When $R$ contains a field, the result was proved in \cite{ELS} and \cite{HH0}.

      \subsubsection{}   An efficient and unified way of dealing with questions related to the homological conjectures in char. $p$ is provided by ``tight closure theory", which has some flavor of almost ring theory. 
      Using Th. \ref{T2} and Th. \ref{T3} above, Heitmann and Ma give evidence that the ``extended plus closure" introduced in \cite{He0} is a good analog of tight closure theory in mixed characteristic \cite{HeMa2}.

     \end{sloppypar}

    \end{document}